\documentclass[12pt]{article}
\usepackage{amsmath,amsfonts,amsthm,amssymb,upgreek}
\usepackage{hyperref}
\usepackage{xcolor}

\newtheorem{thm}{Theorem}

\newtheorem{lem}[thm]{Lemma}

\newtheorem{ass}{Assumption}
\newtheorem{de}{Definition}
\newtheorem{ex}{Example}
\newtheorem{rem}{Remark}

\newcommand{\be}{\begin{equation*}}
\newcommand{\ben}{\begin{equation}}
\newcommand{\ee}{\end{equation*}}
\newcommand{\een}{\end{equation}}

\newcommand{\cc}{{\mathrm c}}

\newcommand{\bbr}{\mathbb{R}}

\newcommand{\gu}{{\mathfrak u}}

\newcommand{\EXP}{\operatorname{\textsf{\upshape E}}}
\newcommand{\PR}{\operatorname{\textsf{\upshape P}}}

\newcommand{\acal}{\mathcal{A}}

\newcommand{\ecal}{\mathcal{E}}
\newcommand{\fcal}{\mathcal{F}}

\newcommand{\di}{\mathrm{d}}

\begin{document}

\title{A risk-sensitive ergodic singular stochastic
control problem}

\author{{\sc Justin Gwee}\footnote{Department of Mathematics,
London School of Economics, Houghton Street, London
WC2A 2AE, UK, \texttt{J.Gwee@lse.ac.uk}.}
\ and
{\sc Mihail Zervos}\footnote{Department of Mathematics,
London School of Economics, Houghton Street, London
WC2A 2AE, UK, \texttt{mihalis.zervos@gmail.com}.}
}

\maketitle

%===============================================
\begin{abstract}
We consider a two-sided singular stochastic control
problem with a risk-sensitive ergodic criterion.
In particular, we consider a stochastic system whose
uncontrolled dynamics are modelled by a linear diffusion.
The control that can be applied to the system is modelled
by an additive finite variation process.
The objective of the control problem is to minimise a
risk-sensitive long-term average criterion that penalises
deviations of the controlled process from a given
interval, as well as the expenditure of control effort.
The stochastic control problem has been partly
motivated by the problem faced by a central bank
who wish to control the exchange rate between its
domestic currency and a foreign currency so that
this fluctuates within a suitable target zone.
We derive the complete solution to the problem under
general assumptions by deriving a $C^2$ solution to
its HJB equation.
To this end, we use the solutions to a suitable
family of Sturm-Liouville eigenvalue problems.
\end{abstract}

%=========================================
\section{Introduction}
We consider a stochastic dynamical system whose state
process satisfies the SDE
\ben
\di X_t = b(X_t) \, \di t + \di \xi_t + \sigma (X_t) \, \di W_t ,
\quad X_0 = x \in \bbr, \label{SDE}
\een
where $W$ is a standard one-dimensional Brownian motion
and $\xi$ is a controlled c\`{a}gl\`{a}d finite-variation process.
With each controlled process $\xi$, we associate the
risk-sensitive long-term average performance index
\ben
J_x (\theta, \xi) = \limsup _{T \uparrow \infty} \frac{1}{\theta T}
\ln \EXP \bigl[ \exp \bigl( \theta I_T (\xi) \bigr) \bigr] , \label{J}
\een
where $\theta > 0$ is the risk-sensitivity parameter
and
\ben
I_T (\xi) = \int _0^T h(X_t) \, \di t + \int _0^T k_+ (X_t)
\oplus \di \xi _t^+ + \int _0^T k_- (X_t) \ominus \di \xi _t^- .
\label{IT}
\een
Here,
\begin{align}
\int _0^T k_+ (X_t) \oplus \di \xi _t^+ & = \int _0^T
k_+ (X_t) \, \di \xi _t^{\cc+} + \sum _{0 \leq t \leq T} \int
_0^{\Delta \xi _t^+} k_+ (X_t + r) \, \di r
\label{k-dxi+} \\
\text{and} \quad
\int _0^T k_- (X_t) \ominus \di \xi _t^- & = \int _0^T
k_- (X_t) \, \di \xi _t^{\cc-} + \sum _{0 \leq t \leq T} \int
_0^{\Delta \xi _t^-} k_- (X_t - r) \, \di r , \label{k-dxi-}
\end{align}
where $\xi^{\cc+}$, $\xi^{\cc-}$ are the continuous parts
of the increasing processes $\xi^+$, $\xi^-$ providing the
unique decomposition $\xi = \xi ^+ - \xi ^-$ and $|\xi|
= \xi ^+ + \xi ^-$, with $|\xi|$ denoting the total variation
process of $\xi$.
The objective of the resulting ergodic risk-sensitive
singular stochastic control problem is to minimise (\ref{J})
over all admissible controlled processes $\xi$.

This stochastic control problem has been partly
motivated by the problem faced by a central bank
who wish to control the exchange rate between its
domestic currency and a foreign currency so that
this fluctuates within a suitable target zone.
In this context, the state process $X$ models the
log exchange rate's stochastic dynamics, while
the controlled process $\xi$ models the cumulative
effect of the bank's interventions in the FX market
to buy or sell the foreign currency.
Furthermore, the running cost function $h$ penalises
deviations of the log exchange rate from a desired
nominal value, while the functions $k_+$ and $k_-$
model proportional transaction costs resulting from
the bank's interventions.

Similar models, which endogenise an exchange
rate's target zone by formulating its management as
a singular stochastic stochastic control problem,
have been studied by
Jeanblanc-Picqu\'{e}~\cite{JP93},
Mundaca and {\O}ksendal~\cite{MO98},
Cadenillas and Zapatero~\cite{CZ99, CZ00},
Ferrari and Vargiolu~\cite{FV20},
and references therein.
The stochastic control problems solved in these references
involve expected discounted performance criteria.
Discounting is commonly used to estimate the present
value of an asset or to model an economic agent's
impatience.
Since an exchange rate is not an asset and a
central bank can be viewed as an institution as well as
a regulator, a long-term average criterion may be more appropriate
for this kind of applications.

Singular stochastic control problems have been motivated
by several applications in areas including target tracking,
optimal harvesting, optimal investment in the presence of
proportional transaction costs and others.
Singular stochastic control problems with risk-neutral
ergodic criteria have been studied by
Karatzas~\cite{K83},
Menaldi and Robin~\cite{MR84, MR13},
Taksar, Klass and Assaf~\cite{TKA88},
Menaldi, Robin and Taksar~\cite{MRT92},
Weerasinghe~\cite{W02, W07},
Jack and Zervos~\cite{JZ06},
L{\o}kka and Zervos~\cite{LZ11, LZ13},
Hynd~\cite{H12},
Wu and Chen~\cite{WC17},
Hening, Nguyen, Ungureanu and Wong~\cite{HNUW19},
Alvarez and Hening~\cite{AH22},
Kunwai, Xi, Yin and Zhu~\cite{KXYZ22},
Liang, Liu and Zervos~\cite{LLZ25},
listed in rough chronological order, and
references therein.
Cohen, Hening and Sun~\cite{CHS22} have also
solved a stochastic game that arises in the
context of ergodic singular stochastic control
with model ambiguity.
On the other hand, Park~\cite[Chapter~I]{P96} and
Chala~\cite{C21} study singular stochastic control problems
with finite time horizon risk-sensitive criteria.
In the context of this paper, Park~\cite[Chapter~II]{P96}
studies a risk-sensitive singular stochastic control problem
with an ergodic criterion in $\bbr^n$, but with constant
$\sigma$.
In this reference, the existence of a suitable solution
to the problem's HJB equation is established and
a limiting connection with the solution to a certain
deterministic ergodic differential game is
established.
For other ergodic risk-sensitive control problems, see
the recent review paper by Biswas and
Borkar~\cite{BB23}.

We derive the complete solution to the problem that
we consider by deriving a $C^2$ solution to the
problem's HJB equation that determines the
optimal strategy, which reflects the state process
in the endpoints of an interval $[\alpha_\star,
\beta_\star]$.
To this end, we first use a suitable logarithmic transformation
that gives rise to a family of Sturm-Liouville eigenvalue
problems parametrised by their boundary points
$\alpha < \beta$.
We then use the optimality conditions suggested by
the so-called smooth-fit of singular stochastic control
to derive the optimal free-boundary points
$\alpha_\star < \beta_\star$.
Furthermore, we show that the control problem's
optimal growth rate identifies with the maximal
eigenvalue of the corresponding Sturm-Liouville
problem.

%=========================================
\section{Problem formulation}
\label{sec:form}

Fix a filtered probability space $\bigl( \Omega, \fcal,
(\fcal_t), \PR \bigr)$ satisfying the usual conditions
and supporting a standard one-dimensional
$(\fcal_t)$-Brownian motion $W$.
We consider a dynamical system, the uncontrolled stochastic
dynamics of which are modelled by the SDE
\ben
\di \overline{X}_t = b(\overline{X}_t) \, \di t + \sigma
(\overline{X}_t) \, \di W_t , \quad \overline{X}_0 = x
\in \bbr . \label{SDE0}
\een
We make the following assumption, which also ensures
that (\ref{SDE0}) has a unique strong solution up to a
possible explosion time.

\begin{ass} \label{A1} {\rm
The functions $b , \sigma :  \bbr \rightarrow \bbr$
are $C^1$ and there exists a constant $C > 0$
such that
\ben
0 < \sigma^2 (x) < C \quad \text{for all }
x \in \bbr . % \label{sigma-ass}
\een
} \end{ass}

We next consider the stochastic control problem defined
by (\ref{SDE})--(\ref{k-dxi-}).

\begin{de} {\rm
The family of all admissible control strategies $\acal$ is the
set of all finite variation $(\fcal_t)$-adapted process $\xi$
with c\`{a}gl\`{a}d sample paths such that $\xi_0=0$ and
the SDE (\ref{SDE}) has a unique non-explosive strong
solution such that
\ben
\limsup _{T \uparrow \infty} \frac{1}{T} \ln
\EXP \Bigl[ \exp \bigl( p |X_T| \bigr) \Bigr] = 0
\quad \text{for all } p>0 . \label{adm-lim}
\een
} \end{de}

\begin{ex} \label{ex:OU} {\rm
Suppose that
\be
\di \overline{X}_t = \upgamma( \upmu - \overline{X}_t)
\, \di t + \upsigma \, \di W_t , \quad \overline{X}_0
= x \in \bbr ,
\ee
for some constants $\upgamma, \upsigma > 0$
and $\upmu \in \bbr$.
Given any $T>0$, the random variable $\overline{X}_T$
has the normal distribution with mean $m_T$ and variance
$\Sigma_T^2$ given by
\be
m_T = \upmu + (x-\upmu) e^{-\upgamma T}
\quad \text{and} \quad
\Sigma_T^2 = \frac{\upsigma^2}{2\upgamma}
(1 - e^{-2\upgamma T}) .
\ee
In view of this observation and the symmetry of the
normal distribution, we can see that
\begin{align}
\EXP \Bigl[ e^{p |\overline{X}_T|} \Bigr] & \leq
e^{p |m_T|} \EXP \Bigl[ e^{p |\overline{X}_T - m_T|}
\Bigr] \leq 2 e^{p |m_T|} \EXP \Bigl[
e^{p (\overline{X}_T - m_T)} \Bigr] \nonumber \\
& = \frac{2}{\sqrt{2\pi \Sigma_T^2}} \exp \biggl(
p |m_T| + \frac{1}{2} p^2 \Sigma_T^2 \biggr)
\xrightarrow [T \uparrow \infty]{}
\sqrt{\frac{4 \upgamma}{\pi \upsigma^2}}
\exp \biggl( p |\upmu| + \frac{p^2 \upsigma^2}{4\upgamma}
\biggr) \nonumber
\end{align}
for all $p>0$, which implies that the choice
$\xi = 0$ is admissible because $\overline{X}$
satisfies (\ref{adm-lim}).
} \end{ex}

\begin{rem} {\rm
In the previous example, the controlled process
$\xi = 0$ is admissible.
However, this is not necessarily true in other
special cases.
For instance, if $\overline{X}$ is a standard Brownian
motion, as in Example~\ref{ex:BM} below, then
the choice $\xi = 0$ is {\em not\/} admissible.
} \end {rem}

We also make the following assumption.
Its requirements on the functions $H_-$ and
$H_+$ are straighforward adaptations of
Assumption~2.3 in Weerasinghe~\cite{W02}
and Assumptions~2.2.(e),(f) in Jack and
Zervos~\cite{JZ06}, who solve risk-neutral
versions of the problem that we study here.
Furthermore, it is of a similar nature as
Assumption~2.9.(ii) in Ferrari and
Vargiolu~\cite{FV20}, who solve a related
singular stochastic control problem with an
expected discounted criterion.
Indeed, this type of an assumption is essential
for the optimal controlled strategy $\xi^\star$
to reflect the associated state process
$X^\star$ at the endpoints of a finite
interval.

\begin{ass} \label{A2} {\rm
The function $h$ is $C^1$ and positive, while
the functions $k_+$ and $k_-$ are $C^2$
and such that
\ben
0 < k_+ (x) < K \quad \text{and} \quad
0 < k_- (x) < K \quad \text{for all } x \in \bbr ,
\label{k-bound}
\een
for some constant $K>0$.
Furthermore, if we define
\begin{align}
H_- (x,\theta) & = \frac{1}{2} \theta
\sigma^2 (x) k_+^2 (x)  - \frac{1}{2} \sigma^2 (x)
k_+' (x) - b(x) k_+ (x) + h(x) \label{H-} \\
\text{and} \quad
H_+ (x,\theta) & = \frac{1}{2} \theta
\sigma^2 (x) k_-^2 (x) + \frac{1}{2} \sigma^2 (x)
k_-' (x) + b(x) k_- (x) + h(x) , \label{H+}
\end{align}
for $x \in \bbr$, then
\ben
\lim _{x \downarrow -\infty} H_- (x, \theta) =
\lim _{x \uparrow \infty} H_+ (x, \theta) = \infty
\quad \text{for all } \theta > 0 \label{Hpm-lims}
\een
and there exist points $\upalpha _- =
\upalpha _- (\theta) \leq \upalpha _+ (\theta)
= \upalpha _+$ such that
\ben
\text{the function } H_- (\cdot , \theta) \text{ is strictly }
\begin{cases} \text{decreasing and positive in }
]{-\infty}, \upalpha _-[ , \\ \text{negative in } ]\upalpha _-,
\upalpha _+[ , \textrm{ if } \upalpha _- < \upalpha_+ , \\
\textrm{increasing and positive in } ]\upalpha _+, \infty[ ,
\end{cases} \label{h-ass1}
\een
as well as constants $\upbeta _- = \upbeta _-
(\theta) \leq \upbeta _+ (\theta) = \upbeta _+$
such that
\ben
\text{the function } H_+ (\cdot , \theta) \text{ is strictly }
\begin{cases} \text{decreasing and positive in } ]{-\infty},
\upbeta _-[ , \\ \text{negative in } ]\upbeta _-, \upbeta _+[
, \text{ if } \upbeta _- < \upbeta _+ , \\ \textrm{increasing
and positive in } ]\upbeta _+, \infty[ , \end{cases}
\label{h-ass2}
\een
for all $\theta > 0$.
} \end{ass}

\begin{ex} {\rm
Let $\overline{X}$ be the process considered in
Example~\ref{ex:OU}.
Also, let $h(x)=cx^2$ and $k_+ (x) = k_- (x) = K$
for some constants $c , K > 0$.
In this context, the functions $H_-$ and $H_+$
defined by (\ref{H-}) and (\ref{H+}) admit the
expressions
\begin{align}
H_- (x, \theta) & = c \biggl( x +\frac{\upgamma K}{2c}
\biggr)^2 + \frac{1}{2} \theta \upsigma^2 K^2 -
\frac{1}{4c} \upgamma^2 K^2 - \upgamma \upmu K
\nonumber \\
\text{and} \quad
H_+ (x, \theta) & = c \biggl( x - \frac{\upgamma K}{2c}
\biggr)^2 + \frac{1}{2} \theta \upsigma^2 K^2 -
\frac{1}{4c} \upgamma^2 K^2 + \upgamma \upmu K .
\nonumber
\end{align}
If $\frac{1}{2} \theta \upsigma^2 K^2 -
\frac{1}{4c} \upgamma^2 K^2 - \upgamma \upmu K
\geq 0$, then $\upalpha_- (\theta) = \upalpha_+
(\theta) = - \frac{1}{2c} \upgamma K$, otherwise
\be
\upalpha_\pm (\theta) = - \frac{\upgamma K}{2c}
\pm \sqrt{- \frac{1}{c} \biggl( \frac{1}{2} \theta
\upsigma^2 K^2 - \frac{1}{4c} \upgamma^2 K^2 -
\upgamma \upmu K \biggr)} .
\ee
Similarly, if $\frac{1}{2} \theta \upsigma^2 K^2 -
\frac{1}{4c} \upgamma^2 K^2 + \upgamma \upmu K
\geq 0$, then $\upbeta_- (\theta) = \upbeta_+
(\theta) = \frac{1}{2c} \upgamma K$, otherwise,
\be
\upbeta_\pm (\theta) = - \frac{\upgamma K}{2c}
\pm \sqrt{- \frac{1}{c} \biggl( \frac{1}{2} \theta
\upsigma^2 K^2 - \frac{1}{4c} \upgamma^2 K^2 +
\upgamma \upmu K \biggr)} .
\ee
In particular, the conditions required by
Assumption~\ref{A2} are all satisfied.
} \end{ex}

\begin{ex} \label{ex:BM} {\rm
Suppose that $\overline{X} = x + \upsigma W$.
Also, let $h(x)=cx^2$ and $k_+ (x) = k_- (x) = K$
for some constants $c , K > 0$.
In this case, the functions $H_-$ and $H_+$
defined by (\ref{H-}) and (\ref{H+}) are given by
\be
H_- (x, \theta) = H_+ (x, \theta) =
cx^2 + \frac{1}{2} \theta \upsigma^2 K^2 .
\ee
and the conditions required by Assumption~\ref{A2}
are all satisfied.
} \end{ex}

%=========================================
\section{The control problem's HJB equation and its
associated Sturm-Liouville eigenvalue problem}
\label{sec:HJB}

Fix any value for the risk-sensitivity parameter
$\theta > 0$.
We will solve the control problem that we consider
by constructing a function $w(\cdot, \theta)$ and
finding a constant $\lambda (\theta)$ such that
$w(\cdot, \theta)$ is $C^2$ and the HJB equation
\begin{align}
\min \biggl\{ \frac{1}{2} \sigma^2 (x) w_{xx} (x, \theta)
+ \frac{1}{2} \theta \bigl( \sigma (x) w_x (x, \theta)
\bigr)^2 + b(x) w_x (x, \theta) + h(x) - \lambda , &
\nonumber \\
k_+ (x) + w_x (x, \theta) , \ k_- (x) - w_x (x, \theta)
& \biggr\} = 0 \label{HJB}
\end{align}
holds true for all $x \in \bbr$.
Given such a solution to this HJB equation,
\be
\inf _{\xi \in \acal} J_x (\theta, \xi) = \lambda (\theta)
\quad \text{for all } x \in \bbr ,
\ee
where $J_x$ is defined by (\ref{J}).
Furthermore, an optimal strategy can be characterised
as follows.
The controller should wait and take no action for as long
as the state process $X$ takes values in the set where
$- k_+ (x) < w_x (x, \theta) < k_- (x)$.
Otherwise, the controller should take minimal action
to keep the state process $X$ outside the interior of
the set in which $w_x (x, \theta) = - k_+ (x)$ or
$w_x (x, \theta) = k_- (x)$ at all times.

We will prove that the optimal control strategy is
characterised by two points $\alpha = \alpha (\theta)
< \beta (\theta) =\beta$ and takes the following form.
If the initial state $x$ is strictly greater than $\beta$
(resp., strictly less than $\alpha$), then it is optimal to
push the state process in an impulsive way down to level
$\beta$ (resp., up to level $\alpha$).
Beyond such a possible initial jump, it is optimal to take
minimal action to keep the state process $X$ inside the
set $[\alpha, \beta]$ at all times, which amounts to
reflecting $X$ in $\beta$ in the negative direction and in
$\alpha$ in the positive direction.
In view of the discussion in the previous paragraph,
the optimality of such a strategy is associated with a
solution $\bigl( w (\cdot, \theta) , \lambda (\theta)
\bigr)$ to the HJB equation (\ref{HJB}) such that
\begin{gather}
w_x (x, \theta) = - k_+ (x) , \quad \text{for } x \in
\mbox{} ]{-\infty}, \alpha] , \label{HJB3} \\
\frac{1}{2} \sigma^2 (x) w_{xx} (x, \theta) +
\frac{1}{2} \theta \bigl( \sigma (x) w_x (x, \theta)
\bigr)^2 + b(x) w_x (x, \theta) + h(x) - \lambda (\theta)
= 0 , \quad \text{for } x \in \mbox{} ]\alpha, \beta[
, \label{HJB1} \\
\text{and} \quad
w_x (x, \theta) = k_- (x) , \quad \text{for } x \in
[\beta, \infty[ . \label{HJB2}
\end{gather}
To determine the points $\alpha < \beta$, we consider the
so-called ``smooth pasting'' condition of singular stochastic
control, which suggests that $w (\cdot, \theta)$ should
be $C^2$, in particular, at the free-boundary points
$\alpha$ and $\beta$.
This condition gives rise to the equations
\begin{gather}
\lim _{x \downarrow \alpha} w_x (x, \theta) =
- k_+ (\alpha) , \quad \lim _{x \downarrow \alpha}
w_{xx} (x, \theta) =  - k_+' (\alpha) , \label{a-eq1} \\
\lim _{x \uparrow \beta} w_x (x, \theta) = k_- (\beta)
\quad \text{and} \quad \lim _{x \uparrow \beta}
w_{xx} (x, \theta) = k_-' (\beta) . \label{b-eq1}
\end{gather}
In view of (\ref{HJB1}), these free-boundary equations
can be satisfied if and only if
\ben
H_- (\alpha, \theta) = \lambda (\theta) = H_+ (\beta, \theta)
, \label{lambda2}
\een
where the functions $H_-$ and $H_+$ are defined by
(\ref{H-}) and (\ref{H+}).

The ODE (\ref{HJB1}) is a Riccati equation.
If we write
\ben
w_x (x, \theta) = \frac{u_x (x, \theta)}{\theta u(x, \theta)}
, \quad \text{for } x \in \mbox{} ]\alpha, \beta[ ,
\label{wandu}
\een
for some function  $u (\cdot, \theta) > 0$,
then $w (\cdot, \theta)$ is a solution to the ODE
(\ref{HJB1}) if and only if $u (\cdot, \theta)$ is a
solution to the second order linear ODE
\be
\frac{1}{2} \sigma^2 (x) u_{xx} (x, \theta) + b(x)
u_x (x, \theta) + \theta \bigl( h(x) - \lambda (\theta)
\bigr) u(x, \theta) = 0 ,
\ee
which is equivalent to
\ben
\frac{\partial}{\partial x} \bigl( q(x) u_x (x, \theta)
\bigr) + \frac{2 \theta}{\sigma^2 (x)} \bigl( h(x)
- \lambda (\theta) \bigr) q(x) u(x, \theta) = 0 ,
\label{ODEu2}
\een
where
\ben
q(x) = \exp \biggl( \int _0^x \frac{2b(y)}{\sigma^2 (y)}
\, \di y \biggr) . \label{q}
\een
In view of this transformation and the boundary conditions
(\ref{a-eq1}) and (\ref{b-eq1}), we are faced with the
regular Sturm-Liouville eigenvalue problem defined
by the ODE (\ref{ODEu2}) and the boundary conditions
\ben
\theta k_+ (\alpha) u(\alpha, \theta) + u_x (\alpha
, \theta) = 0 \quad \text{and} \quad
\theta k_- (\beta) u(\beta, \theta) - u_x (\beta ,
\theta) = 0 . \label{SV-BC}
\een
This problem has infinitely many simple real eigenvalues
\be
\lambda_0 (\theta) > \lambda_1 (\theta) > \cdots
> \lambda_n (\theta) > \cdots \quad \text{such that }
\lim _{n \uparrow \infty} \lambda_n (\theta)
= -\infty
\ee
and no other eigenvalues, while the eigenfunction
$\gu^{(n)} (\cdot, \theta)$ corresponding to $\lambda_n
(\theta)$ has exactly $n$ zeros in the interval
$]\alpha, \beta[$ (e.g., see
Walter~\cite[Theorem~VI.27.II]{W}).
Furthermore, the eigenvalues are related to their
corresponding eigenfunctions by means of the Rayleigh
quotient
\begin{align}
\lambda_n (\theta) = \Biggl( & q(\beta) \gu^{(n)}
(\beta, \theta) \gu_x^{(n)} (\beta, \theta) - q(\alpha)
\gu^{(n)} (\alpha, \theta) \gu_x^{(n)} (\alpha, \theta)
\nonumber \\
& + \int _\alpha^\beta q(y) \biggl( \frac{2 \theta h(y)}
{\sigma^2 (y)} \bigl( \gu^{(n)} (y, \theta) \bigr) ^2
- \bigl( \gu_x^{(n)} (y, \theta) \bigr) ^2 \biggr)
\, \di y \Biggr) \nonumber \\
& \times \biggl( \int _\alpha^\beta
\frac{2\theta q(y)}{\sigma^2 (y)} \bigl( \gu^{(n)}
(y, \theta) \bigr) ^2 \, \di y \biggr) ^{-1} .
\label{Rayleigh}
\end{align}

The eigenfunction $\gu^{(0)} (\cdot, \theta)$ is the
only one that has no zeros in $]\alpha, \beta[$.
The function $w_x (\cdot, \theta)$ given by
(\ref{wandu}) is therefore clearly well-defined only
for $u (\cdot, \theta) = \gu^{(0)} (\cdot, \theta)$.
In view of this observation, we consider the maximal
eigenvalue $\lambda_0 (\theta)$ and its corresponding
eigenfunction $\gu^{(0)} (\cdot, \theta)$ in what follows.
We also write $\uplambda (\alpha, \beta, \theta)$
and $\phi _{\alpha, \beta, \theta}$ instead of
$\lambda_0 (\theta)$ and $\gu^{(0)} (\cdot, \theta)$
to stress their dependence on the free-boundary
points $\alpha$ and $\beta$, as well as on the risk
sensitivity parameter $\theta$.
Furthermore, we assume that $\phi
_{\alpha, \beta, \theta}$ has been normalised by
a multiplicative constant, so that
\ben
\int _\alpha^\beta \frac{2\theta q(y)}{\sigma^2 (y)}
\phi _{\alpha, \beta, \theta}^2 (y) \, \di y = 1 ,
\label{phi-norm}
\een
and we note that the boundary conditions (\ref{SV-BC})
and the expression (\ref{Rayleigh}) imply that
\begin{align}
\uplambda (\alpha, \beta, \theta) = \mbox{} &
\theta \Bigl( q(\alpha) k_+ (\alpha) \phi
_{\alpha,\beta, \theta}^2 (\alpha) + q(\beta) k_-
(\beta) \phi _{\alpha,\beta, \theta}^2 (\beta)
\Bigr) \nonumber \\
& + \int _\alpha^\beta q(y) \biggl (\frac{2\theta h(y)}
{\sigma^2 (y)} \phi _{\alpha,\beta, \theta}^2 (y)
- \bigl( \phi _{\alpha,\beta, \theta}' (y) \bigr)^2
\biggr) \, \di y . \label{uplambda}
\end{align}

\begin{lem} \label{uplambda-props}
In the presence of Assumptions~\ref{A1} and~\ref{A2},
the function $\uplambda$ defined by (\ref{uplambda})
for $\alpha < \beta$ and $\theta > 0$ is $C^{1,1,1}$,
\begin{gather}
\uplambda _\alpha (\alpha, \beta, \theta) =
\frac{2 \theta q(\alpha)}{\sigma^2 (\alpha)}
\phi _{\alpha,\beta, \theta}^2 (\alpha) \bigl( \uplambda
(\alpha, \beta, \theta) - H_- (\alpha, \theta) \bigr) ,
\label{uplambda_a} \\
\uplambda _\beta (\alpha, \beta, \theta) = -
\frac{2 \theta q(\beta)}{\sigma^2 (\beta)}
\phi _{\alpha, \beta, \theta}^2 (\beta) \bigl( \uplambda
(\alpha, \beta, \theta) - H_+ (\beta, \theta) \bigr)
\label{uplambda_b} \\
\text{and} \quad
\uplambda _\theta (\alpha, \beta, \theta) =
\frac{1}{\theta} \int _\alpha^\beta q(y) \bigl(
\phi _{\alpha, \beta, \theta}' (y) \bigr)^2 \, \di y
> 0 , \label{uplambda_theta}
\end{gather}
where the functions $H_-$ and $H_+$ are defined by
(\ref{H-}) and (\ref{H+})).
\end{lem}
\noindent
{\bf Proof.}
%========
We prove these identities using a technique inspired
by Kong and Zettl~\cite{KZ}.
To establish (\ref{uplambda_a}), we fix any $\theta > 0$
and we drop it from the notation of the functions $H_\pm$,
$\uplambda$ and $\phi$.
Given any $\varepsilon > 0$, we use integration by parts
and the ODE (\ref{ODEu2}) to calculate
\begin{align}
q(\beta) \Bigl( \phi _{\alpha, \beta} (\beta) & \phi
_{\alpha + \varepsilon, \beta}' (\beta) - \phi _{\alpha, \beta}'
(\beta) \phi _{\alpha + \varepsilon, \beta} (\beta) \Bigr)
\nonumber \\
& - q(\alpha {+} \varepsilon) \Bigl( \phi _{\alpha, \beta}
(\alpha {+} \varepsilon) \phi _{\alpha + \varepsilon, \beta}'
(\alpha {+} \varepsilon) - \phi _{\alpha, \beta}' (\alpha
+ \varepsilon) \phi _{\alpha {+} \varepsilon, \beta}
(\alpha {+} \varepsilon) \Bigr) \nonumber \\
= \mbox{} & \int _{\alpha + \varepsilon}^\beta \Bigl( \phi
_{\alpha, \beta} (y) \bigl( q \phi _{\alpha + \varepsilon, \beta}'
\bigr)' (y) - \phi _{\alpha + \varepsilon, \beta} (y) \bigl(
q \phi  _{\alpha, \beta}' \bigr)' (y) \Bigr) \, \di y
\nonumber \\
= \mbox{} & \bigl( \uplambda (\alpha {+} \varepsilon, \beta)
- \uplambda (\alpha, \beta) \bigr) \int _{\alpha + \varepsilon}
^\beta \frac{2\theta q(y)}{\sigma^2 (y)} \phi_{\alpha,\beta}
(y) \phi _{\alpha + \varepsilon, \beta} (y) \, \di y .
\nonumber 
\end{align}
In view of the boundary conditions (\ref{SV-BC}), these
identities imply that
\begin{align}
\bigl( \uplambda (\alpha {+} \varepsilon, \beta)
- \uplambda & (\alpha, \beta) \bigr) \int _{\alpha + \varepsilon}
^\beta \frac{2\theta q(y)}{\sigma^2 (y)} \phi_{\alpha,\beta}
(y) \phi _{\alpha + \varepsilon, \beta} (y) \, \di y
\nonumber \\
& = q(\alpha {+} \varepsilon) \Bigl( \theta k_+ (\alpha
{+} \varepsilon) \phi _{\alpha, \beta} (\alpha {+} \varepsilon)
+ \phi _{\alpha, \beta}' (\alpha {+} \varepsilon) \Bigr)
\phi _{\alpha + \varepsilon, \beta} (\alpha {+} \varepsilon)
. \nonumber
\end{align}
Using the ODE (\ref{ODEu2}) and
the boundary conditions (\ref{SV-BC}) once more,
we obtain
\begin{align}
q(\alpha {+} \varepsilon) & \Bigl( \theta k_+ (\alpha {+}
\varepsilon) \phi _{\alpha, \beta} (\alpha {+} \varepsilon)
+ \phi _{\alpha, \beta}' (\alpha {+} \varepsilon) \Bigr)
\nonumber \\
& = \int _\alpha^{\alpha + \varepsilon} \Bigl( \theta
\bigl( k_+ q \phi _{\alpha, \beta} \bigr)' (y) + \bigl( q \phi
_{\alpha, \beta}' \bigr)' (y) \Bigr) \, \di y \nonumber \\
& = \int _\alpha^{\alpha + \varepsilon}
\frac{2 \theta q(y)}{\sigma^2 (y)} \phi _{\alpha,\beta} (y)
\biggl( \uplambda (\alpha, \beta) + \frac{1}{2}
\sigma^2 (y) k_+ (y) \frac{\phi _{\alpha,\beta}' (y)}
{\phi _{\alpha,\beta} (y)} \nonumber \\
& \hspace{45mm}
+ \frac{1}{2} \sigma^2 (y) k_+ (y) + b(y) k_+ (y)
- h(y) \biggr) \, \di y . \nonumber
\end{align}
It follows that
\begin{align}
& \bigl( \uplambda (\alpha {+} \varepsilon, \beta)
- \uplambda (\alpha, \beta) \bigr) \int
_{\alpha + \varepsilon}^\beta \frac{2\theta q(y)}
{\sigma^2 (y)} \phi_{\alpha,\beta} (y) \phi
_{\alpha + \varepsilon, \beta} (y) \, \di y
\nonumber \\
& = \phi _{\alpha {+} \varepsilon, \beta}
(\alpha + \varepsilon) \int _\alpha^{\alpha + \varepsilon}
\frac{2 \theta q(y)}{\sigma^2 (y)} \phi _{\alpha,\beta} (y)
\biggl( \uplambda (\alpha, \beta) + \frac{1}{2}
\sigma^2 (y) k_+ (y) \biggl( \frac{\phi _{\alpha,\beta}' (y)}
{\phi _{\alpha,\beta} (y)} + \theta k_+ (y) \biggr)
- H_- (y) \biggr) \, \di y . \nonumber
\end{align}
Dividing by $\varepsilon$ and passing to the limit as
$\varepsilon \downarrow 0$ using (\ref{SV-BC}),
as well as (\ref{phi-norm}), we can see that the
right-hand derivative $\uplambda _{\alpha+}
(\alpha, \beta)$ exists and is equal to the expression
on the right-hand side of (\ref{uplambda_a}).

Replacing $\alpha$ and $\alpha+\varepsilon$ by
$\alpha-\varepsilon$ and $\alpha$, respectively, in the
analysis above, we can see that the left-hand derivative
$\uplambda _{\alpha-} (\alpha, \beta)$ also exists and is
equal to $\uplambda _{\alpha+} (\alpha, \beta)$.

The proof of (\ref{uplambda_b}) follows the same arguments.

To prove (\ref{uplambda_theta}), we fix any
$\alpha < \beta$ and we write $\uplambda (\theta)$
and $\phi _\theta$ in place of $\uplambda (\alpha,
\beta, \theta)$ and $\phi _{\alpha, \beta, \theta}$.
Given any $\varepsilon \neq 0$ small, we use the boundary
conditions (\ref{SV-BC}), integration by parts
and the ODE (\ref{ODEu2}) to calculate
\begin{align}
\varepsilon \Bigl( q(\beta) k_- (\beta) & \phi _\theta
(\beta) \phi _{\theta + \varepsilon} (\beta) + q(\alpha)
k_+ (\alpha) \phi _\theta (\alpha) \phi
_{\theta + \varepsilon} (\alpha) \Bigr) \nonumber \\
= \mbox{} & q(\beta) \Bigl( \phi _\theta (\beta)
\phi _{\theta + \varepsilon}' (\beta) - \phi _\theta'
(\beta) \phi _{\theta + \varepsilon} (\beta) \Bigr)
- q(\alpha) \Bigl( \phi _\theta (\alpha) \phi
_{\theta + \varepsilon}' (\alpha) - \phi _\theta'
(\alpha) \phi _{\theta + \varepsilon} (\alpha) \Bigr)
\nonumber \\
= \mbox{} & \int _\alpha^\beta \Bigl(
\phi _\theta (y) \bigl( q \phi _{\theta + \varepsilon}'
\bigr)' (y) - \phi _{\theta + \varepsilon} (y) \bigl(
q \phi  _\theta' \bigr)' (y) \Bigr) \, \di y
\nonumber \\
= \mbox{} & \Bigl( (\theta {+} \varepsilon) \uplambda
(\theta {+} \varepsilon) - \theta \uplambda (\theta)
\Bigr) \int _\alpha^\beta \frac{2q(y)}{\sigma^2 (y)}
\phi _\theta (y) \phi _{\theta + \varepsilon} (y)
\, \di y \nonumber \\
& - \varepsilon \int _\alpha^\beta
\frac{2 q(y) h(y)}{\sigma^2 (y)} \phi _\theta (y)
\phi _{\theta + \varepsilon} (y) \, \di y . \nonumber 
\end{align}
Dividing by $\varepsilon$ and passing to the
limit as $\varepsilon \downarrow 0$, we obtain
\begin{align}
\biggl( \uplambda' (\theta) + \frac{1}{\theta}
\uplambda (\theta) \biggr) \int _\alpha^\beta
& \frac{2\theta q(y)}{\sigma^2 (y)} \phi _\theta^2
(y) \, \di y \nonumber \\
& = q(\beta) k_- (\beta) \phi _\theta^2 (\beta)
+ q(\alpha) k_+ (\alpha) \phi _\theta^2 (\alpha) +
\int _\alpha^\beta \frac{2q(y) h(y)}{\sigma^2 (y)}
\phi _\theta^2 (y) \, \di y . \nonumber
\end{align}
Combining this result with (\ref{phi-norm})
and (\ref{uplambda}), we obtain
(\ref{uplambda_theta}).
\mbox{}\hfill$\Box$
\bigskip

We prove the following result here, rather than in the
context of Theorem~\ref{thm:main}, because we will
need the strict positivity of the function $\uplambda$
to derive the solution to the HJB equation
(\ref{HJB}) that identifies the optimal strategy.

\begin{lem} \label{lem:J-lam}
Suppose that Assumptions~\ref{A1} and~\ref{A2}
hold true.
The function $\uplambda$ defined by (\ref{uplambda})
is such that, given any points $\alpha < \beta$ in $\bbr$,
\ben
\uplambda (\alpha, \beta, \theta) = J_x \bigl( \theta,
\xi^{\alpha, \beta} \bigr) > 0 \quad \text{for all }
x \in \bbr , \label{SL-lambda-J}
\een
where $J_x$ is defined by (\ref{J}) and $\xi^{\alpha, \beta}
\in \acal$ is the controlled process that,
beyond an initial jump $\Delta \xi_0^{\alpha, \beta}
= (\alpha-x)^+ - (x-\beta)^+$, is continuous and
reflects the corresponding state process
$X^{\alpha, \beta}$ in $\alpha$ in the positive direction
and in $\beta$ in the negative direction.
\end{lem}
\noindent
{\bf Proof.}
%========
Formally, the controlled process $\xi^{\alpha, \beta}$
and the corresponding solution $X^{\alpha, \beta}$
to the SDE (\ref{SDE}) are characterised by the
requirements that
\begin{gather}
X_T^{\alpha, \beta} \in [\alpha, \beta] , \quad
\xi _T^{\alpha, \beta, +} - (\alpha-x)^+ =
\int _0^T {\bf 1} _{\{ X_t^{\alpha, \beta} = \alpha \}}
\, \di \xi _t^{\alpha, \beta, \cc +} \label{X-xi-ab1} \\
\text{and} \quad
\xi _T^{\alpha, \beta, -} - (x-\beta)^+ = \int _0^T
{\bf 1} _{\{ X_t^{\alpha, \beta} = \beta \}}
\, \di \xi _t^{\alpha, \beta, \cc -} \label{X-xi-ab2}
\end{gather}
for all $T>0$.
Such processes indeed exists (e.g., see
Tanaka~\cite[Theorem~4.1]{T79}).
In particular, $\xi^{\alpha, \beta}$ belongs
to $\acal$ because $X_t^{\alpha, \beta} \in
[\alpha, \beta]$ for all $t > 0$.

To establish (\ref{SL-lambda-J}), we first consider any
$C^1$ function $w: \bbr \rightarrow \bbr$ that is
piece-wise $C^2$ and any admissible control strategy
$\xi \in \acal$.
Using It\^{o}-Tanaka's formula for general semimartingales
and the identities $\Delta X_t = X_{t+} - X_t = \Delta
\xi _t$, we obtain
\begin{align}
w(X_{T+}) = \mbox{} & w(x) + \int _0^T \biggl( \frac{1}{2}
\sigma^2 (X_t) w'' (X_t) + b(X_t) w' (X_t) \biggr) \, \di t
+ \int _{[0,T]} w' (X_t) \, \di \xi_t \nonumber \\
& + \sum _{0 \leq t \leq T} \bigl( w(X_{t+}) - w(X_t)
- w' (X_t) \Delta X_t \bigr) + M_T \nonumber \\
= \mbox{} & w(x) + \int _0^T \biggl( \frac{1}{2}
\sigma^2 (X_t) w'' (X_t) + b(X_t) w' (X_t) \biggr) \, \di t
+ \int _0^T w' (X_t) \, \di \xi _t^{\cc+} \nonumber \\
& - \int _0^T w' (X_t) \, \di \xi_t^{\cc-} + \sum
_{0 \leq t \leq T} \bigl( w(X_{t+}) - w(X_t) \bigr) + M_T
, \nonumber
\end{align}
where
\ben
M_T = \int _0^T \sigma (X_t) w'(X_t) \, \di W_t .
\label{M}
\een
Combining this expression with the identity
\be
w(X_{t+}) - w(X_t) = \int _0^{\Delta\xi _t^+}
w'(X_t+r) \, \di r - \int _0^{\Delta \xi _t^-}
w' (X_t-r) \, \di r ,
\ee
we can see that
\begin{align}
\sum _{0 \leq t \leq T} & \bigl( w(X_{t+}) - w(X_t)
\bigr) + \sum _{0 \leq t \leq T} \int _0^{\Delta \xi _t^+}
k_+ (X_t + r) \, \di r + \sum _{0 \leq t \leq T}
\int _0^{\Delta \xi _t^-} k_- (X_t - r) \, \di r
\nonumber \\
& = \sum _{0 \leq t \leq T} \int _0^{\Delta \xi _t^+} \bigl(
k_+ (X_t + r) + w' (X_t+r) \bigr) \, \di r + \sum _{0 \leq t \leq T}
\int _0^{\Delta \xi _t^-} \bigl( k_- (X_t - r) - w' (X_t-r) \bigr)
\, \di r . \nonumber
\end{align}
Recalling the definitions (\ref{k-dxi+}) and (\ref{k-dxi-}),
we obtain
\begin{align}
\int _0^T & h(X_t) \, \di t + \int _0^T k_+ (X_t) \oplus
\di \xi _t^+ + \int _0^T k_- (X_t) \ominus \di \xi _t^-
+ w(X_{T+}) \nonumber \\
= \mbox{} & \uplambda (\alpha, \beta) T + w(x)
- \frac{1}{2} \theta \langle M \rangle_T
+ M_T \nonumber \\
& + \int _0^T \biggl( \frac{1}{2} \sigma^2 (X_t) w''(X_t)
+ b(X_t) w'(X_t) + \frac{1}{2} \theta \bigl( \sigma
(X_t) w'(X_t)\bigr)^2 + h(X_t) - \uplambda
(\alpha, \beta) \biggr) \, \di t \nonumber \\
& + \int _0^T \bigl( k_+ (X_t) + w'(X_t) \bigr) \,\di
\xi _t^{\cc+} + \int _0^T \bigl( k_- (X_t) - w'(X_t)
\bigr) \,\di \xi _t^{\cc-} \label{Ito-lem1} \\
& + \sum _{0 \leq t \leq T} \int _0^{\Delta \xi _t^+}
\bigl( k_+ (X_t+r) + w'(X_t+r) \bigr) \, \di r +
\sum _{0 \leq t \leq T} \int _0^{\Delta \xi _t^-} \bigl(
k_- (X_t-r) - w'(X_t-r) \bigr) \, \di r . \nonumber
\end{align}

Let $w$ be a function whose first derivative is given by
\be
w' (x) = \begin{cases} - k_+ (x) , & \text{if } x \leq
\alpha , \\ \frac{1}{\theta} \frac{\di}{\di x} \ln \bigl( \phi
_{\alpha, \beta, \theta} (x) \bigr) , & \text{if } x \in \mbox{}
]\alpha, \beta[ , \\ k_- (x) , & \text{if } x \geq \beta
. \end{cases}
\ee
Recalling that the eigenfunction $\gu_0 (\cdot, \theta) =
\phi _{\alpha,\beta,\theta}$ and the eigenvalue $\lambda_0
= \uplambda (\alpha, \beta, \theta)$ provide a solution to
the Sturm-Liouville eigenvalue problem defined by the
ODE (\ref{ODEu2}) with boundary conditions (\ref{SV-BC}),
we can see that $w$ satisfies (\ref{HJB1}).
Furthermore, in view of (\ref{wandu}), (\ref{SV-BC}) and
the in-between arguments, we can see that $w$ is
$C^2$ in $\bbr \setminus \{ \alpha, \beta \}$ and
$C^1$ at both of $\alpha$ and $\beta$. 
Combining these observations with (\ref{X-xi-ab1}),
(\ref{X-xi-ab2}) and (\ref{Ito-lem1}), we obtain
\ben
I_T (\xi _t^{\alpha, \beta})
= \uplambda (\alpha, \beta) T + w(x) - w \bigl(
X_{T+}^{\alpha, \beta} \bigr) - \frac{1}{2} \theta
\bigl\langle M^{\alpha, \beta} \bigr\rangle_T +
M_T^{\alpha, \beta} , \label{IT-ab}
\een
where $M^{\alpha, \beta}$ is defined by (\ref{M})
for $X = X^{\alpha, \beta}$.

The assumption  (\ref{k-bound}) and the definition
of $w'$ imply that $w'$ is bounded.
Combining this observation with the assumption
that $\sigma$ is bounded, we can see that
there exists a constant $C_1 > 0$ such that
$\langle M^{\alpha, \beta} \rangle_T
\leq C_1 T$ for all $T>0$.
Therefore, the process defined by
\ben
\ecal _T \bigl( \theta M^{\alpha, \beta} \bigr)
= \exp \biggl( - \frac{1}{2} \theta^2 \bigl\langle
M^{\alpha, \beta} \bigr\rangle _T
+ \theta M_T^{\alpha, \beta} \biggr)
\label{exp-mart}
\een
is a martingale, thanks to Novikov's condition.
In view of this observation and (\ref{IT-ab}),
we obtain
\begin{align}
\frac{1}{\theta T} \ln \EXP \bigl[ \exp \bigl( \theta
I_T (\xi _t^{\alpha, \beta} \bigr) \bigr] 
= \mbox{} & \uplambda (\alpha, \beta) + \frac{w(x)}{T}
\nonumber \\
& + \frac{1}{\theta T} \ln \EXP \Biggl[ \exp \biggl(
\theta \biggl( -w \bigl( X_T^{\alpha, \beta} \bigr)
- \frac{1}{2} \theta \bigl\langle M^{\alpha, \beta}
\bigr\rangle_T + M_T^{\alpha, \beta} \biggr)
\biggr) \Biggr]
\nonumber \\
= \mbox{} & \uplambda (\alpha, \beta) + \frac{w(x)}{T}
+ \frac{1}{\theta T} \ln
\EXP^{\widetilde{\PR}_T^{\alpha,\beta}}
\Bigl[ \exp \Bigl( - \theta w \bigl(X_T^{\alpha, \beta}
\bigr) \Bigr) \Bigr] , \nonumber
\end{align}
where $\widetilde{\PR}_T^{\alpha,\beta}$ is the
probability measure on $(\Omega, \fcal_T)$ with
Radon-Nikodym derivative with respect to $\PR$
given by $\di \widetilde{\PR}_T^{\alpha,\beta}
/ \di {\PR} = \ecal _T \bigl( \theta M^{\alpha, \beta}
\bigr)$.
Finally, using the fact that the process
$w(X^{\alpha, \beta})$ is bounded, we can pass
to the limit as $T \uparrow \infty$ to obtain
the identity
$J_x (\theta, \xi^{\alpha, \beta}) = \uplambda
(\alpha, \beta, \theta)$.
\mbox{}\hfill$\Box$

%=========================================
\section{The solution to the control problem}
\label{sec:sol}

The following result identifies the solution to the
HJB equation (\ref{HJB}) that yields the control
problem's solution.

\begin{thm} \label{HJBsoln}
In the presence of Assumptions~\ref{A1} and~\ref{A2},
the following statements hold true.
\smallskip

\noindent
{\rm{(I)}}
Given any $\theta > 0$, there exists a unique pair
$\bigl( \alpha_\star (\theta) , \beta_\star (\theta) \bigr)$
such that
\begin{gather}
\alpha_\star (\theta) < \upalpha_- (\theta) ,
\quad \upbeta_+ (\theta) < \beta_\star (\theta)
\label{ab*} \\
\text{and} \quad
\lambda_\star (\theta) := \uplambda \bigl(
\alpha_\star (\theta), \beta_\star (\theta), \theta
\bigr) = H_- \bigl( \alpha_\star (\theta), \theta
\bigr) = H_+ \bigl( \beta_\star (\theta), \theta
\bigr) , \label{lam*}
\end{gather}
where the function $\uplambda$ is defined by
(\ref{uplambda}), the points $\upalpha_- (\theta)$,
and $\upbeta_+ (\theta)$ are as in (\ref{h-ass1})
and (\ref{h-ass2}), while the functions $H_-$
and $H_+$ are defined by (\ref{H-}) and (\ref{H+}).
\smallskip

\noindent
{\rm{(II)}}
The function $w(\cdot, \theta)$ that is defined by
\ben
w_x (x, \theta) = \begin{cases}
- k_+ (x), & \text{if } x \leq \alpha_\star (\theta) , \\
\frac{1}{\theta} \frac{\di}{\di x} \ln \bigl( \phi
_{\alpha_\star (\theta) , \beta_\star (\theta) , \theta}
(x) \bigr) , & \text{if } x \in \mbox{} \bigl]
\alpha_\star (\theta) , \beta_\star (\theta) \bigr[ , \\
k_- (x) , & \text{if } x \geq \beta_\star (\theta) ,
\end{cases} \label{wdef}
\een
modulo an additive constant, is $C^2$.
Furthermore, this function and $\lambda_\star
(\theta)$ provide a solution to the HJB equation
(\ref{HJB}).
\end{thm}
\noindent
{\bf Proof.}
%=======
Throughout the proof, we fix any $\theta > 0$ and
we drop it from the notation of the functions
$H_\pm$, $\uplambda$, $\phi$, $\alpha_\star$,
$\beta_\star$ and $w$.
\smallskip

{\em Proof of (I).\/}
The conditions (\ref{Hpm-lims})--(\ref{h-ass2}) in
Assumption~\ref{A2} ensure the existence of a unique
function $\Gamma: [\upbeta_+ , \infty[ \mbox{}
\rightarrow \mbox{} ]{-\infty} , \upalpha_-]$ such that
$H_+ (\beta) = H_- \bigl( \Gamma (\beta) \bigr)$ for
all $\beta \geq \upbeta_+$.
In particular, $\Gamma (\upbeta_+) = \upalpha_-$.
The $C^1$ continuity of the functions $b$, $\sigma$ and
$h$, together with the $C^2$ continuity of the functions
$k_+$ and $k_-$, implies that the both of the functions
$H_-$ and $H_+$ defined by (\ref{H-}) and (\ref{H+})
are $C^1$ (see Assumptions~\ref{A1} and~\ref{A2}).
Consequently, the restriction of $\Gamma$ to
the interval $]\upbeta_+ , \infty[$ is also $C^1$.
In light of these observations, if the equation
\ben
\Lambda (\beta) := \uplambda \bigl( \Gamma
(\beta), \beta \bigr) = H_+ (\beta) \label{beta*-eqn}
\een
admits a unique solution $\beta_\star > \upbeta_+$, then
part~(I) of the theorem holds with $\alpha_\star
= \Gamma (\beta_\star)$.

To show that the equation (\ref{beta*-eqn}) has a unique
solution $\beta_\star > \upbeta_+$, we first use
(\ref{uplambda_a}) and (\ref{uplambda_b}) in
Lemma~\ref{uplambda-props}, as well as the
identity $H_+ (\beta) = H_- \bigl( \Gamma (\beta) \bigr)$,
to calculate
\begin{align}
\frac{\di}{\di \beta} \bigl( \Lambda & (\beta) - H_+ (\beta)
\bigr) \nonumber \\
&= \uplambda_\alpha \bigl( \Gamma (\beta) , \beta
\bigr) \Gamma' (\beta) + \uplambda_\beta \bigl(
\Gamma (\beta) , \beta) - H_+' (\beta)
\nonumber \\
& = 2 \theta \Biggl(
\frac{q \bigl( \Gamma (\beta) \bigr) \phi
_{\Gamma (\beta) , \beta}^2 \bigl( \Gamma (\beta) \bigr)}
{\sigma^2 \bigl( \Gamma (\beta) \bigr)} \Gamma' (\beta)
- \frac{q(\beta) \phi _{\Gamma (\beta) , \beta}^2 (\beta)}
{\sigma^2 (\beta)} \Biggr)
\bigl( \Lambda (\beta) - H_+ (\beta) \bigr) - H_+' (\beta)
\nonumber \\
& =: \varrho (\beta) \bigl( \Lambda(\beta) - H_+ (\beta)
\bigr) - H_+' (\beta) , \quad \text{for } \beta > \upbeta_+
. \label{Lambda-ODE}
\end{align}
The solution to this first-order ODE is such that
\begin{align}
I(\beta) \bigl( \Lambda (\beta) - H_+ (\beta) \bigr)
& = \Lambda (\upbeta_+) - H_+ (\upbeta_+) - \int
_{\upbeta_+}^\beta I(u) H_+' (u) \, \di u \nonumber \\
& = \uplambda (\upalpha_-, \upbeta_+) - \int
_{\upbeta_+}^\beta I(u) H_+' (u) \, \di u
=: F(\beta) , \nonumber
\end{align}
where $I(\beta) = \exp \bigl( - \int _{\upbeta_+}^\beta
\varrho (u) \, \di u \bigr)$.
The second equality here follows from the fact that
$\Gamma (\upbeta_+) = \upalpha_-$ and the assumption
that $H_+ (\upbeta_+) = 0$.
It follows that equation (\ref{beta*-eqn}) is equivalent
to the equation
\ben
F(\beta) = 0 . \label{F-eqn}
\een

In view of the inequalities
\be
F' (\beta) = -I(\beta) H_+' (\beta) < 0 \text{ for all }
\beta > \upbeta_+ \quad \text{and} \quad
F(\upbeta_+) = \uplambda (\upalpha_-, \upbeta_+)
\stackrel{(\ref{SL-lambda-J})}{>} 0 ,
\ee
we can see that equation (\ref{F-eqn}) has a unique
solution $\beta_\star > \upbeta_+$ if and only if
$\lim _{\beta \uparrow \infty} F(\beta) < 0$.
To see that this inequality is indeed true, we argue by
contradiction.
To this end, we assume that
$\lim _{\beta \uparrow \infty} F(\beta) \geq 0$,
which can be true only if
\ben
\Lambda (\beta) - H_+ (\beta) =
\frac{F(\beta)}{I(\beta)} > 0 \quad \text{for all }
\beta > \upbeta_+ \label{Lambda-h}
\een
because $F' (\beta) < 0$ and $I(\beta) > 0$
for all $\beta > \upbeta_+$.
In view of the inequalities $\Gamma' < 0$
and $q > 0$, we can see that the function
$\varrho$ introduced in (\ref{Lambda-ODE})
is such that $\varrho (\beta) < 0$ for all $\beta
> \upbeta_+$.
In view of this inequality, the contradiction
hypothesis (\ref{Lambda-h}) and the identity
\be
\Lambda' (\beta) = \varrho (\beta) \bigl( \Lambda
(\beta) - H_+ (\beta) \bigr) ,
\ee
which follows from (\ref{Lambda-ODE}), we can
see that $\Lambda' (\beta) < 0$ for all $\beta >
\upbeta_+$.
However, this conclusion and (\ref{Hpm-lims})
imply that 
\be
\lim _{\beta \uparrow \infty} \bigl( \Lambda (\beta)
- H_+ (\beta) \bigr) \leq \Lambda(\upbeta_+)
- \lim _{\beta \uparrow \infty} H_+ (\beta) = -\infty ,
\ee
which contradicts (\ref{Lambda-h}).
Thus, we have proved that equation (\ref{F-eqn}),
which is equivalent to equation (\ref{beta*-eqn}),
has a unique solution $\beta_\star > \upbeta_+$
and we have established part~(I) of the theorem.
\smallskip

{\em Proof of (II).\/}
By construction, we will prove that the function $w$
given by \eqref{wdef} is a $C^2$ solution to the
HJB equation (\ref{HJB}) if we show that
\begin{gather}
\frac{1}{2} \theta \sigma^2 (x) k_+^2 (x) 
- \frac{1}{2} \sigma^2 (x) k_+' (x) - b(x) k_+ (x) + h(x)
- \lambda_\star \geq 0 \quad \text{for all }
x < \alpha_\star , \label{HJB-ineq11} \\
\frac{1}{2} \theta \sigma^2 (x) k_-^2 (x) +
\frac{1}{2} \sigma^2 (x) k_-' (x) + b(x) k_- (x) + h(x)
- \lambda_\star \geq 0 \quad \text{for all }
x > \beta_\star \label{HJB-ineq12} \\
\text{and} \quad - k_+ (x) \leq w'(x) \leq k_- (x) \quad
\text{for all } x \in \mbox{} ]\alpha_\star, \beta_\star[ .
\label{HJB-ineq2}
\end{gather}
The inequalities (\ref{HJB-ineq11}) and
(\ref{HJB-ineq12}) follow immediately from
(\ref{h-ass1}) and (\ref{h-ass2}) in Assumption~\ref{A2}
once we observe that
\begin{gather}
\frac{1}{2} \theta \sigma^2 (x) k_+^2 (x) 
- \frac{1}{2} \sigma^2 (x) k_+' (x) - b(x) k_+ (x)
+ h(x) - \lambda_\star = H_- (x) - H_- (\alpha_\star)
\quad \text{for all } x < \alpha_\star
\nonumber \\
\intertext{and}
\frac{1}{2} \theta \sigma^2 (x) k_-^2 (x) +
\frac{1}{2} \sigma^2 (x) k_-' (x) + b(x) k_- (x)
+ h(x) - \lambda_\star = H_+ (x) - H_+
(\beta_\star) \quad \text{for all } x > \beta_\star
, \nonumber
\end{gather}
where we have used the definitions (\ref{H-}) and
(\ref{H+}) of the functions $H_-$ and $H_+$, as well
as part~(I) of the theorem.

To establish (\ref{HJB-ineq2}), we first note that
the $C^1$ continuity of the functions $b$, $\sigma$
and $h$ implies that the restriction of $w$ in
$]\alpha_\star, \beta_\star[$ is $C^3$.
In particular, we note that differentiation of the ODE
(\ref{HJB1}) that $w$ satisfies in $]\alpha_\star,
\beta_\star[$ implies that
\begin{align}
\frac{1}{2} \sigma^2 (x) w''' (x) + \bigl( b(x)
+ \sigma (x) \sigma' (x) + \theta \sigma^2 (x) w'(x)
\bigr) w'' (x) & \nonumber \\
\mbox{} + \theta \sigma (x) \sigma' (x) \bigl(
w'(x) \bigr)^2 + b'(x) w'(x) + h'(x) & = 0 . \nonumber
\end{align}
In view of this calculation, the inequalities
(\ref{ab*}), the assumptions (\ref{h-ass1}),
(\ref{h-ass2}) and the free-boundary equations
(\ref{a-eq1}), (\ref{b-eq1}), we can see that
\begin{align}
\lim _{x \downarrow \alpha_\star} \bigl( w''' (x)
+ k_+'' (x) \bigr) & = - \frac{2}{\sigma^2 (\alpha_\star)}
H_-' (\alpha_\star) > 0 \nonumber \\
\text{and} \quad
\lim _{x \uparrow \beta_\star} \bigl( w''' (x)
- k_-'' (x) \bigr) & = - \frac{2}{\sigma^2 (\beta_\star)}
H_+' (\beta_\star) < 0 . \nonumber
\end{align}
It follows that there exists $\varepsilon > 0$ such
that
\begin{align}
w''(x) + k_+' (x) & > 0 \quad \text{for all } x \in
\mbox{} ]\alpha_\star, \alpha_\star + \varepsilon[
\label{w''(a*close)} \\
\text{and} \quad
w''(x) - k_-' (x) & > 0 \quad \text{for all } x \in
\mbox{} ]\beta_\star - \varepsilon , \beta_\star[ .
\label{w''(b*close)} 
\end{align}

We next argue by contradiction, we assume
that there exist $x \in \mbox{} ]\alpha_\star,
\beta_\star[$ such that $w'(x) > k_- (x)$ and
we define
\begin{align}
\alpha_\star < \underline{\gamma} & := \min \bigl\{
x \in \mbox{} ]\alpha_\star, \beta_\star[ \mbox{}
\mid \ w'(x) = k_- (x) \bigr\} \nonumber \\
& < \max \bigl\{ x \in \mbox{} ]\alpha_\star,
\beta_\star[ \mbox{} \mid \ w'(x) = k_- (x)
\bigr\} =: \overline{\gamma} < \beta_\star ,
\label{HJB-ineq2-c}
\end{align}
where the inequalities follow once we combine the
boundary conditions $w'(\alpha_\star) = - k_+
(\alpha_\star)$ and $w'(\beta_\star) = k_-
(\beta_\star)$ with (\ref{w''(b*close)}).
Combining the definitions of the points
$\underline{\gamma}$ and $\overline{\gamma}$
in (\ref{HJB-ineq2-c}) with the boundary conditions
$w'(\alpha_\star) = - k_+ (\alpha_\star)$ and
$w'(\beta_\star) = k_- (\beta_\star)$, we can
see that
\be
w'' (\underline{\gamma}) - k_-' (\underline{\gamma})
\geq 0 \quad \text{and} \quad
w'' (\overline{\gamma}) - k_-' (\overline{\gamma})
\leq 0 .
\ee
On the other hand, using the ODE (\ref{HJB1}),
the definitions (\ref{H-}) and (\ref{H+}) of the functions
$H_-$ and $H_+$, and part~(I) of the theorem,
we obtain
\begin{align}
w'' (\underline{\gamma}) - k_-' (\underline{\gamma})
& = \frac{2}{\sigma^2 (\underline{\gamma})}
\bigl( H_+ (\beta_\star) - H_+ (\underline{\gamma})
\bigr) \nonumber \\
\text{and} \quad
w'' (\overline{\gamma}) - k_-' (\overline{\gamma})
& = \frac{2}{\sigma^2 (\overline{\gamma})}
\bigl( H_+ (\beta_\star) - H_+ (\overline{\gamma})
\bigr) . \nonumber
\end{align}
However, these inequalities and expressions
associated with the points $\underline{\gamma}
< \overline{\gamma} < \beta_\star$ contradict
(\ref{h-ass2}) in Assumption~\ref{A2},
and the right-hand side of (\ref{HJB-ineq2})
follows.

Finally, we can show that the left-hand side
of (\ref{HJB-ineq2}) holds true using
(\ref{w''(a*close)}) and a contradiction argument
similar to the one based on (\ref{HJB-ineq2-c}).
\mbox{}\hfill$\Box$
\bigskip

We can now prove the main result of the paper.

\begin{thm} \label{thm:main}
Suppose that Assumptions~\ref{A1} and~\ref{A2}
hold true.
If $\bigl( \alpha_\star (\theta) , \beta_\star (\theta)
\bigr)$ and $\lambda_\star (\theta)$ are as in
Theorem~\ref{HJBsoln}, then, given any $x \in \bbr$,
\ben
\inf _{\xi \in \acal} J_x (\theta, \xi) = J_x \bigl( \theta,
\xi^\star \bigr) = \lambda_\star (\theta) > 0 ,
\label{J-thm}
\een
where the controlled process $\xi^\star \in \acal$ is
continuous beyond an initial jump $\Delta \xi_0^\star
= \bigl( \alpha_\star (\theta) -x \bigr)^+ -
\bigl( x - \beta_\star (\theta) \bigr)^+$ and reflects the
corresponding state process $X^\star$ in
$\alpha_\star (\theta)$ in the positive direction
and in $\beta_\star (\theta)$ in the negative direction.
\end{thm}
\noindent
{\bf Proof.} 
%=======
Fix any $x \in \bbr$, $\theta > 0$ and $\xi \in \acal$.
Also, consider the solution to the HJB equation
(\ref{HJB}) presented by Theorem~\ref{HJBsoln}.
Given $\varepsilon \in \mbox{} ]0 , \theta[$, the
expression (\ref{Ito-lem1}) in the proof of
Lemma~\ref{lem:J-lam} with $w (\cdot, \theta {-}
\varepsilon)$ in place of $w$ implies that
\be
I_T (\xi) - \lambda_\star (\theta {-} \varepsilon) T
+ w(X_{T+} , \theta {-} \varepsilon) \geq w(x
, \theta {-} \varepsilon) - \frac{1}{2} (\theta {-}
\varepsilon) \bigl\langle M \bigr\rangle_T
+ M_T ,
\ee
where $M$ is defined by (\ref{M}) for $w' =
w_x (\cdot , \theta  {-} \varepsilon)$.
The exponential local martingale
$\ecal \bigl( (\theta {-} \varepsilon) M \bigr)$ that
is defined by (\ref{exp-mart}) with $M$ in place
of $M^{\alpha, \beta}$ is a martingale because
$\sigma$ and $w_x (\cdot, \theta)$ are both
bounded (see also the discussion above
(\ref{exp-mart})).
Therefore,
\be
\EXP \Biggl[ \exp \biggl( - \frac{1}{2} (\theta {-}
\varepsilon) ^2 \bigl\langle M \bigr\rangle_T
+ (\theta {-} \varepsilon) M_T \biggr) \Biggr] = 1 .
\ee
In view of these observations and H\"{o}lder's
inequality, we can see that
\begin{align}
\exp \bigl( (\theta {-} \varepsilon) w(x , \theta {-}
\varepsilon) \bigr) & \leq \EXP \bigl[ \exp \bigl(
(\theta {-} \varepsilon) I_T (\xi) - (\theta {-} \varepsilon)
\lambda_\star (\theta {-} \varepsilon) T + (\theta {-}
\varepsilon) w(X_T , \theta {-} \varepsilon) \bigr)
\bigr] \nonumber \\
& \leq \Bigl( \EXP \bigl[ \exp \bigl( \theta I_T (\xi)
- \theta \lambda_\star (\theta {-} \varepsilon) T
\bigr) \bigr] \Bigr) ^{\frac{\theta {-} \varepsilon}{\theta}}
\Bigl( \EXP \bigl[ \exp \bigl( \varepsilon^{-1} \theta
(\theta {-} \varepsilon) w(X_T , \theta {-} \varepsilon)
\bigr) \bigr] \Bigr) ^{\frac{\varepsilon}{\theta}}
. \nonumber
\end{align}
It follows that
\begin{align}
\frac{w(x , \theta {-} \varepsilon)}{T} & \leq
\frac{1}{\theta T} \ln \EXP \bigl[ \exp \bigl( \theta I_T
(\xi) \bigr) \bigr] - \lambda_\star (\theta {-} \varepsilon)
+ \frac{\varepsilon}{\theta (\theta {-} \varepsilon) T}
\ln \EXP \biggl[ \exp \biggl(
\frac{\theta (\theta {-} \varepsilon) K}{\varepsilon}
|X_T| \biggr) \biggr] .
\nonumber
\end{align}
Recalling the admissibility condition (\ref{adm-lim})
and passing to the limit as $T \uparrow \infty$ in
this inequality, we obtain $J_x (\theta, \xi) \geq
\lambda_\star (\theta {-} \varepsilon)$.
The inequality $J_x (\theta, \xi) \geq \lambda_\star
(\theta)$ follows by passing to the limit as
$\varepsilon \downarrow 0$ because
$\lambda_\star$ is continuous.

Finally, the identity $J_x \bigl( \theta,
\xi^\star \bigr) = \lambda_\star (\theta)$ and
the optimality of the controlled process
$\xi^\star$ follow from Lemma~\ref{lem:J-lam}.
\mbox{}\hfill$\Box$
\bigskip

We conclude the paper with the following result
on the dependence of the control problem’s
solution on the risk-sensitivity parameter
$\theta$.

\begin{lem}
In the presence of Assumptions~\ref{A1}
and~\ref{A2}, the following statements hold true.
\smallskip

\noindent {\rm (I)}
The optimal growth rate $\lambda_\star$
is such that
\be
\lambda_\star' (\theta) > 0
\quad \text{and} \quad
\lim _{\theta \downarrow 0}
\lambda_\star' (\theta) = \infty .
\ee

\noindent {\rm (II)}
If $\sigma$ is constant and $k_+ (x)
= k_- (x) = K$ for some constant $K > 0$, then
the free-boundaries $\alpha_\star < \beta_\star$
are such that
\be
\alpha _\star' (\theta) < 0
\quad \text{and} \quad
\beta _\star' (\theta) > 0 .
\ee
\end{lem}
\noindent
{\bf Proof.}
%=======
Differentiating the identities (\ref{lam*}), we
obtain
\ben
\lambda_\star' (\theta) =
\frac{\di \uplambda \bigl( \alpha _\star (\theta),
\beta _\star (\theta), \theta \bigr)}{\di \theta} =
\frac{\di H_- \bigl( \alpha _\star (\theta) , \theta \bigr)}
{\di \theta} =
\frac{\di H_+ \bigl( \beta _\star (\theta) , \theta \bigr)}
{\di \theta} . \label{lamH-derivs}
\een
Using the partial derivatives given by
(\ref{uplambda_a}), (\ref{uplambda_b}) and
(\ref{uplambda_theta}), as well as the identities
(\ref{lam*}), we can see that
\be
\frac{\di \uplambda \bigl( \alpha _\star (\theta),
\beta _\star (\theta), \theta \bigr)}{\di \theta} =
\uplambda _\theta \bigl( \alpha _\star (\theta),
\beta _\star (\theta), \theta \bigr) .
\ee
Part~(I) of the lemma follows from this identity
and the expression (\ref{uplambda_theta})
for $\uplambda _\theta$.

In the context of part~(II) of the lemma, we use the
definitions (\ref{H-}) and (\ref{H+}) of the functions
$H_-$ and $H_+$ to calculate
\begin{align}
\frac{\di H_- \bigl( \alpha _\star (\theta) , \theta \bigr)}
{\di \theta} & =
\frac{\partial H_- \bigl( \alpha _\star (\theta) , \theta \bigr)}
{\partial x} \alpha _\star' (\theta) + \frac{1}{2} K^2
\sigma^2 \nonumber \\
\text{and} \quad
\frac{\di H_+ \bigl( \beta _\star (\theta) , \theta \bigr)}
{\di \theta} & =
\frac{\partial H_+ \bigl( \beta _\star (\theta) , \theta \bigr)}
{\partial x} \beta _\star' (\theta) + \frac{1}{2} K^2
\sigma^2 . \nonumber
\end{align}
These identities and the last equality in
(\ref{lamH-derivs}) imply that
\be
\frac{\partial H_- \bigl( \alpha _\star (\theta) , \theta \bigr)}
{\partial x} \frac{\di \alpha _\star (\theta)}{\di \theta}
= \frac{\partial H_+ \bigl( \beta _\star (\theta) , \theta \bigr)}
{\partial x} \frac{\di \beta _\star (\theta)}{\di \theta} .
\ee
The claims of part~(II) of the lemma follows from this result,
(\ref{h-ass1}) and (\ref{h-ass2}) in Assumption~\ref{A2}
and the inequalities (\ref{ab*}).
\mbox{}\hfill$\Box$

%===============================================
\section*{Acknowledgment}

We are grateful to two anonymous referees and
the Associate Editor whose thoughtful and constructive
comments significantly improved the paper.

%=========================================

\end{document}